\documentclass[10pt, reqno]{amsart}
\usepackage{amssymb}

\newtheorem{thm1}{Theorem}[section]
\newtheorem{lem1}[thm1]{Lemma}
\newtheorem{rem1}[thm1]{Remark}
\newtheorem{def1}[thm1]{Definition}

\newtheorem{prop1}[thm1]{Proposition}
\newtheorem{ex1}[thm1]{Example}
\DeclareMathOperator{\supp}{supp}

\begin{document}

\title[Arithmetical rank of toric ideals associated to graphs]
{Arithmetical rank of toric ideals associated to graphs}

\author[Anargyros Katsabekis]{ Anargyros Katsabekis}
\address { Department of Mathematics, University of the Aegean,
83200 Karlovassi, Samos, GREECE } \email{katsabek@aegean.gr}


\begin{abstract}
\par Let $I_{G} \subset K[x_{1},\ldots,x_{m}]$ be the toric ideal associated to a finite graph $G$. In this paper we study the binomial arithmetical rank and the $G$-homogeneous arithmetical rank of $I_G$ in 2 cases:
\begin{enumerate} \item $G$ is bipartite, \item $I_G$ is generated
by quadratic binomials. \end{enumerate} In both cases we prove
that the binomial arithmetical rank and the $G$-arithmetical rank
coincide with the minimal number of generators of $I_G$.
\end{abstract}
 \maketitle

\section{Introduction} Let $G$ be a finite, connected and undirected graph having no loops and no multiple
edges on the vertex set $V(G)=\{v_{1},\ldots,v_{n}\}$, and let
$E(G)=\{e_1,\ldots,e_m\}$ be the set of edges of $G$. The {\em
incidence matrix} of $G$ is the $n \times m$ matrix
$M_G=(a_{i,j})$
defined by $$a_{i,j}= \left \{\begin{array} {lll} 1, & \textrm{if} \ v_{i} \ \textrm{is one of the vertices in} \ e_{j}\\
 0, & \textrm{otherwise}. \end{array} \right.$$ Let $K$ be an algebraically closed field and let $A_G=\{{\bf a}_1,\ldots,{\bf a}_m\}$ be the set of vectors in $\mathbb{Z}^n$, where ${\bf a}_i=(a_{1,i},\ldots,a_{n,i})$ for $1 \leq i \leq m$. The {\em toric ideal} $I_{G}$ associated to $G$ is the kernel of the $K$-algebra
homomorphism
$$ \phi: K[x_1,\dots, x_m]\rightarrow
K[t_1,\dots,t_{n}]$$ given by
$$\phi(x_i)=t_{1}^{a_{1,i}} \cdots t_{n}^{a_{n,i}} \qquad
\mbox{for all }i = 1,\dots ,m.$$ The ideal $I_{G}$ is prime and
therefore $rad(I_{G})=I_{G}$. The {\em toric variety}
$\mathbb{V}(I_{G})$ associated to $G$ is the set
$$\{(u_1,\ldots,u_m) \in K^m|F(u_1,\ldots,u_m)=0, \forall F \in
I_{G}\}$$of zeroes of $I_G$. For every graph $G$ the variety
$\mathbb{V}(I_{G})$ is an {\em extremal toric variety}, i.e. the
vector configuration $A_G$ is extremal, see Remark 2.1.
\par The polynomial ring
$K[x_1,\ldots,x_m]$ has a natural $G$-graded structure given by
setting $\deg_{G}(x_i)={\bf a}_i$ for $i=1,\ldots,m$. For ${\bf
u}=(u_1,\ldots,u_m) \in \mathbb{N}^m$, we define the $G$-{\em
degree} of the monomial ${\bf x}^{{\bf u}}:=x_{1}^{u_1} \cdots
x_{m}^{u_m}$ to be \[ \deg_{G}({\bf x}^{{\bf u}}):=u_1{\bf
a}_1+\cdots+u_m{\bf a}_m \in \mathbb{N}A_{G},\] where
$\mathbb{N}A_{G}$ is the semigroup generated by $A_G$. Remark that
$\mathbb{N}A_{G}$ is {\em pointed}, i.e. zero is the only
invertible element. A polynomial $F\in K[x_1,\dots ,x_m]$ is
called $G$-{\em homogeneous} if the monomials in non-zero terms of
$F$ have the same $G$-degree. An ideal $I$ is $G$-homogeneous if
it is generated by $G$-homogeneous polynomials.\\ The toric ideal
$I_G$ is generated by all the binomials ${\bf x}^{{\bf u} _+}-
{\bf x}^{{\bf u} _-}$ such that $\deg_{G}({\bf x}^{{\bf u}
_+})=\deg_{G}({\bf x}^{{\bf u} _-})$, where ${\bf u} _+\in
\mathbb{N}^m$ and ${\bf u} _-\in \mathbb{N}^m$ denote the positive
and negative part of ${\bf u}={\bf u}_+ -  {\bf u}_-$,
respectively (see \cite{St}).
\par A basic problem in Commutative Algebra asks to compute the
smallest integer $s$ for which there exist polynomials
$F_1,\ldots, F_s$ in $I_G$ such that $I_G=rad(F_1,\ldots,F_s)$.
This integer is called the {\em arithmetical rank} of $I_G$ and
will be denoted by ${\rm ara}(I_G)$. An usual approach to this
problem is to restrict to a certain class of polynomials and ask
how many polynomials from this class can generate the toric ideal
up to radical. Restricting the polynomials to the class of
binomials we arrive at the notion of the {\em binomial
arithmetical rank} of $I_G$, denoted by ${\rm bar}(I_G)$. Also, if
all of the polynomials $F_1,\ldots ,F_s$ satisfying
$I_G=rad(F_1,\ldots,F_s)$ are $G$-homogeneous, the smallest
integer $s$ is called {\em the} $G$-{\em homogeneous arithmetical
rank} of $I_G$ and will be denoted by ${\rm ara}_{G}(I_G)$. From
the definitions, the generalized Krull's principal ideal theorem
and the graded version of Nakayama's Lemma we deduce the following
inequalities for a toric ideal $I_{G}$:
$${\rm ht}(I_{G})\leq {\rm ara}(I_{G})\leq {\rm ara}_{G}
(I_{G})\leq {\rm bar}(I_{G}) \leq \mu(I_{G}).$$ Here ${\rm
ht}(I_G)$ denotes the height and $\mu(I_G)$ denotes the minimal
number of generators of $I_G$. When ${\rm ht}(I_G)=\mu(I_{G})$ the
ideal $I_G$ is called a {\em complete intersection}. \par A case
of particular interest is when ${\rm bar}(I_{G})={\rm ht}(I_{G})$.
When $K$ is a field of characteristic zero this is equivalent to
say that $I_{G}$ is complete intersection, see \cite{BMT}.
Complete intersection bipartite graphs have been characterized in
\cite{GRV}, \cite{Katzman}. In most cases, when $G$ is bipartite,
the equality ${\rm bar}(I_{G})={\rm ht}(I_{G})$ does not hold. In
section 3 we prove that
${\rm bar}(I_{G})={\rm ara}_{G}(I_{G})=\mu(I_{G})$, for any bipartite graph $G$. In addition we show that the equality ${\rm bar}(I_{G})={\rm ara}_{G}(I_{G})=\mu(I_{G})$ also holds, for any graph $G$ such that the toric ideal $I_{G}$ is generated by quadratic binomials.\\

\section{Basics on toric ideals associated to graphs}

\par Let $G$ be a graph. A {\em walk} of length $q$ of $G$ is a finite sequence of the form $$\Gamma=(\{v_{1},v_{2}\}, \{v_{2},v_{3}\},\ldots,\{v_{q-1},v_{q}\},\{v_{q},v_{q+1}\});$$this walk is {\em closed} if $v_{q+1}=v_{1}$. An {\em even} closed walk is a closed walk of even
length. A {\em cycle} of $G$ is a closed walk
$$\Gamma=(\{v_{1},v_{2}\},
\{v_{2},v_{3}\},\ldots,\{v_{q},v_{1}\})$$ with $v_{i} \neq v_{j}$
for all $1 \leq i<j \leq q$. Notice that if
$e_{i}=\{v_{i_k},v_{i_l}\}$ is an edge of $G$, then
$\phi(x_{i})=t_{i_k}t_{i_l}$. Given an even closed walk
$\Gamma=(e_{i_1},\ldots,e_{i_{2q}})$ of $G$ with each $e_{k} \in
E(G)$, we have that
$$\phi(\prod_{k=1}^{q}x_{i_{2k-1}})=\phi(\prod_{k=1}^{q}x_{i_{2k}})$$
and therefore the binomial
$$f_{\Gamma}:=\prod_{k=1}^{q}x_{i_{2k-1}}-\prod_{k=1}^{q}x_{i_{2k}}$$
belongs to $I_{G}$. Remark that if $\Gamma$ is an even cycle of
$G$, then the monomials $M=\prod_{k=1}^{q}x_{i_{2k-1}}$ and
$N=\prod_{k=1}^{q}x_{i_{2k}}$ are squarefree. From Proposition
3.1 in \cite{Vil} we have that
$$I_{G}=(\{f_{\Gamma}| \Gamma \ \textrm{is an even closed walk of}
\ G\}).$$
\par Let ${\bf u}=(u_{1},\ldots,u_{m}) \in \mathbb{Z}^m$ be a
vector, then the {\em support} of ${\bf u}$, denoted by
$\supp({\bf u})$, is the set $\{i \in \{1,\ldots,m\} | u_i \neq
0\}$. For a monomial ${\bf x}^{{\bf u}}$ we define $\supp({\bf
x}^{{\bf u}}):=\supp({\bf u})$. A non-zero vector ${\bf
u}=(u_1,\ldots,u_m)\in {\ker}_{\mathbb{Z}}(M_G)$ is called a {\em
circuit} of $A_G$ if its support is minimal with respect to
inclusion and all the coordinates of ${\bf u}$ are relatively
prime, where ${\ker}_{\mathbb{Z}}(M_G)=\{{\bf v} \in
\mathbb{Z}^{m} | M_{G}{\bf v}^{t}={\bf 0}^{t}\}$. The binomial
${\bf x}^{{\bf u}_+}- {\bf x}^{{\bf u} _-} \in I_G$ associated to
a vector ${\bf u} \in {\ker}_{\mathbb{Z}}(M_G)$ is called also
circuit. A binomial $B={\bf x}^{{\bf u}_+}- {\bf x}^{{\bf u} _-}
\in I_G$ is called {\em primitive} if there exists no other
binomial ${\bf x}^{{\bf v}_+}- {\bf x}^{{\bf v} _-} \in I_G$ such
that ${\bf x}^{{\bf v}_+}$ divides ${\bf x}^{{\bf u}_+}$ and ${\bf
x}^{{\bf v}_-}$ divides ${\bf x}^{{\bf u}_-}$. For a circuit
$B={\bf x}^{{\bf u}_+}- {\bf x}^{{\bf u} _-} \in I_{G}$ we have,
from Corollary 8.1.4 in \cite{V}, that $B=f_{\Gamma}$ for an even
closed walk $\Gamma$ of $G$, since every circuit is also
primitive.\\ If $G$ is a bipartite graph, then $G$ has no odd
cycles, so, from Proposition 4.2 in \cite{Vil}, a binomial
$f_{\Gamma}$, where $\Gamma$ is an even closed walk of $G$, is a
circuit if and only if $\Gamma$ is an even cycle.
\par For the rest of this section we recall some fundamental material from \cite{KT}.\\ We shall denote by
$\mathcal{C}_G$ the set of circuits of $A_G$. Let
$$\mathcal{C}:=\{E \subset \{1,\ldots,m\} \ | \ \supp({\bf u}_+)=E \ {\rm or} \ \supp({\bf u}_-)=E \
{\rm where} \ {\bf u}\in \mathcal{C}_G\}$$ and let
$\mathcal{C}_{min}$ be the set of minimal elements of $\mathcal{C}$.\\
To every toric ideal $I_{G}$ we associate the rational polyhedral
cone $$\sigma=pos_{\mathbb{Q}}(A_G):=\{\lambda_1{\bf
a}_1+\cdots+\lambda_m {\bf a}_m \ | \ \lambda_i \in
\mathbb{Q}_{\geq 0}\}.$$ A {\em face} $\mathcal{F}$ of $\sigma$ is
any set of the form
$$\mathcal{F}= \sigma\cap \{{\bf x}\in \mathbb{Q}^n: {\bf cx} = 0\}, $$ where
${\bf c}\in \mathbb{Q}^n$ and ${\bf cx} \geq 0$ for all ${\bf
x}\in \sigma$. Given an edge $e_{i}=\{v_{i_1},v_{i_2}\}$ of $G$,
we have that $pos_{\mathbb{Q}}({\bf a}_i)$ is a face of $\sigma$
with defining vector ${\bf c}=(c_{1},\ldots,c_{n}) \in \mathbb{Z}^n$ having coordinates $$c_j= \left \{\begin{array} {lll} 0, & \textrm{if} \ j=i_{1},i_{2}\\
 1, & \textrm{otherwise}. \end{array} \right.$$
\begin{rem1} {\rm For every graph $G$ the vector configuration $A_G$ is extremal,
i.e. for any $B\subsetneqq A_G$ we have $pos_{\mathbb{Q}}(B)
\subsetneqq pos_{\mathbb{Q}}(A_G)$.\\ To see this consider a set
$B=\{{\bf a}_{i_1},\ldots,{\bf a}_{i_k}\} \subsetneqq A_G$ and
assume that the vector ${\bf a}_j$ is not in $B$. Let ${\bf c}$ be
the defining vector of the face $pos_{\mathbb{Q}}({\bf a}_j)$. If
${\bf a}_j$ belongs to $pos_{\mathbb{Q}}(B)$, then ${\bf
a}_j=\lambda_{1}{\bf a}_{i_1}+\cdots+\lambda_{k}{\bf a}_{i_k}$,
where $\lambda_{1},\ldots,\lambda_{k}$ are nonnegative rationals
and there is at least one $\lambda_{r}$ different from zero. Thus
${\bf c}{\bf a}_j=\lambda_{1}({\bf c}{\bf
a}_{i_1})+\cdots+\lambda_{k}({\bf c}{\bf a}_{i_k})$ and therefore
$0=\lambda_{1}({\bf c}{\bf a}_{i_1})+\cdots+\lambda_{k}({\bf
c}{\bf a}_{i_k})$. But $$\lambda_{1}({\bf c}{\bf
a}_{i_1})+\cdots+\lambda_{k}({\bf c}{\bf a}_{i_k})>0,$$ a
contradiction. Consequently ${\bf a}_j$ does not belong to
$pos_{\mathbb{Q}}(B)$, so $pos_{\mathbb{Q}}(B)$ is a proper subset
of $pos_{\mathbb{Q}}(A_G)$.}
\end{rem1}

For a subset $E$ of $\{1,\dots ,m\}$ we denote by $\sigma _E$ the
subcone $pos_{\mathbb{Q} }({\bf a}_i |i \in E)$ of $\sigma $. We
adopt the convention that $\sigma _{\emptyset }=\{{\bf 0}\}$. The
{\em relative interior} of $\sigma _E$, denoted by
$relint_{\mathbb{Q}}(\sigma _E)$, is the set of all
 strictly positive rational linear combinations of ${\bf a}_i$, $i \in E$.
\begin{def1} (\cite{KT}) {\rm We associate to $G$ the simplicial complex $\Delta
_{G}$ with vertices the elements  of $\mathcal{C}_{min}$. Let
$T\subset \mathcal{C}_{min}$ then $T\in \Delta _{G}$ if
$$\cap _{E \in T}relint_{\mathbb{Q}}\left(
\sigma_{E}\right) \neq \emptyset .$$In particular $\{E, E'\} \in
\Delta _{G}$ if and only if there exists a circuit ${\bf u} \in
\mathcal{C}_G$ such that $\supp({\bf u}_+)=E$ and $\supp({\bf
u}_-)=E'$.}
\end{def1}
Let $J$ be a subset of $\Omega:=\{0,1,\dots ,{\rm dim}(\Delta
_{G})\}$. A set $\mathcal{M}=\{T_1,\dots ,T_s\}$ of simplices of
$\Delta _{G}$ is called a $J$-{\em matching} in $\Delta _{G}$ if
$T_k\cap T_l=\emptyset$ for every $1 \leq k, l \leq s$ and ${\rm
dim}(T_k) \in J$ for every $1 \leq k \leq s$; see also Defnition
2.1 in \cite{KT}. Let $\supp(\mathcal{M})=\cup^{s}_{i=1}T_i$,
which is a subset of the vertices $\mathcal{C}_{min}$. A
$J$-matching $\mathcal{M}$ in $\Delta _{G}$ is called a {\em
maximal} $J$-matching if $\supp(\mathcal{M})$ has the maximum
possible cardinality among all $J$-matchings.\\Given a maximal
$J$-matching $\mathcal{M}=\{T_1,\dots ,T_s\}$ in $\Delta _{G}$, we
shall denote by $card(\mathcal{M})$ the cardinality $s$ of the set
$\mathcal{M}$. In addition by $\delta(\Delta _{G})_{J}$ we denote
the minimum of the set
$$ \{card(\mathcal{M})| \mathcal{M} \ \textrm{is a maximal} \
J-\textrm{matching} \ \textrm{in} \ \Delta _{G}\}.$$ It follows,
from the definitions, that if $\Delta
_{G}=\bigcup_{i=1}^{t}{\Delta _{G}^{i}}$ then
$$\delta(\Delta _{G})_J=\sum_{i=1}^{t}\delta({\Delta
_{G}^{i}})_J,$$ where $\Delta _{G}^{i}$ are the connected
components of $\Delta _{G}$.

\begin{ex1} {\rm Consider the complete graph $\mathcal{K}_4$ on the vertex set $\{v_{1},\ldots,v_{4}\}$. We consider one variable $x_{ij}$, $1 \leq i<j \leq 4$, for each edge
$\{v_{i},v_{j}\}$ of $\mathcal{K}_4$ and form the polynomial ring
$K[x_{12}, x_{13}, x_{14}, x_{23}, x_{24}, x_{34}]$. From
Proposition 4.2 in \cite{Vil} we have that the toric ideal
$I_{\mathcal{K}_4}$ has 3 circuits, namely
$f_{\Gamma_{1}}=x_{12}x_{34}-x_{14}x_{23}$,
$f_{\Gamma_{2}}=x_{12}x_{34}-x_{13}x_{24}$ and
$f_{\Gamma_{3}}=x_{13}x_{24}-x_{14}x_{23}$, corresponding to the 3
even cycles $\Gamma_{1}$, $\Gamma_{2}$ and $\Gamma_{3}$,
respectively, of $\mathcal{K}_{4}$. In fact $I_{\mathcal{K}_4}$ is
minimally generated by two of the above binomials, so it is
complete intersection of height 2. The simplicial complex
$\Delta_{\mathcal{K}_4}$ has three vertices, namely
$E_{1}=\{12,34\}$, $E_{2}=\{14,23\}$ and $E_{3}=\{13,24\}$. It
consists of all subsets of the set $\{E_{1}, E_{2}, E_{3}\}$.
There are four maximal $\{0,1\}$-matchings in
$\Delta_{\mathcal{K}_4}$, namely $\{\{E_{1},E_{2}\}, \{E_{3}\}\}$,
$\{\{E_{1},E_{3}\}, \{E_{2}\}\}$, $\{\{E_{2},E_{3}\}, \{E_{1}\}\}$
and $\{\{E_1\},\{E_2\},\{E_3\}\}$. We have that $\delta(\Delta
_{\mathcal{K}_4})_{\{0,1\}}=2$ which is attained for example by
the maximal $\{0,1\}$-matching $\{\{E_{1},E_{2}\}, \{E_{3}\}\}$.
In addition $\delta(\Delta _{\mathcal{K}_4})_{\{0,1,2\}}=1$ which
is attained by the maximal $\{0,1,2\}$-matching
$\{\{E_{1},E_{2},E_{3}\}\}$.}
\end{ex1}

Using the fact that $A_G$ is an extremal vector configuration and
also two results from \cite{KT}, namely Theorem 4.6 and Theorem
3.5, we get the following Theorem:
\begin{thm1} For a toric ideal $I_{G}$ we have $\delta(\Delta _{G})_{\{0,1\}} \leq {\rm
bar}(I_G)$ and $\delta(\Delta _{G})_{\Omega} \leq {\rm
ara}_{G}(I_G)$.\\
\end{thm1}

\section{Arithmetical rank}

\par Let $I_G \subset K[x_{1},\ldots,x_{m}]$ be the toric ideal associated to a graph $G$. A binomial $B \in I_{G}$ is called {\em indispensable} if every
system of binomial generators of $I_{G}$ contains $B$ or $-B$,
while a monomial $M$ is called {\em indispensable} if every system
of binomial generators of $I_{G}$ contains a binomial $B$ such
that $M$ is a monomial of $B$. Let $\mathcal{N}_{G}$ be the
monomial ideal generated by all ${\bf x}^{\bf u}$ for which there
exists a nonzero ${\bf x}^{\bf u}-{\bf x}^{\bf w} \in I_G$. From
Proposition 3.1 in \cite{CKT} we have that the set of
indispensable monomials is the unique minimal generating set of
$\mathcal{N}_{G}$. The following lemma will be useful in the proof
of Theorem 3.2 and Proposition 3.4.

\begin{lem1} Assume that either $G$ is a bipartite graph or $I_G$
is generated by quadratic binomials. Let $T_G=\{M_1,\ldots,M_r\}$
be the set of indispensable monomials, then
$\mathcal{C}_{min}=\{\supp(M_1),\ldots,\supp(M_r)\}$.
\end{lem1}

\noindent \textbf{Proof.} Consider first the case that $G$ is a
bipartite graph. From Theorem 3.2 in \cite{Katzman} we have that
$I_G$ is minimally generated by all binomials of the form
$f_{\Gamma}$, where $\Gamma$ is an even cycle of $G$ with no
chord. Combining the above theorem and Theorem 2.3 in \cite{Hi-O}
we obtain that a binomial ${\bf x}^{{\bf u}_+}- {\bf x}^{{\bf u}
_-} \in I_G$ is indispensable if and only if it is of the form
$f_{\Gamma}$, for an even cycle of $G$ with no chord. Notice that
in some cases there are circuits of $I_G$ of the form
$f_{\Gamma}$, for an even cycle $\Gamma$ of $G$ with a chord. If
$B_1,\ldots,B_s$ are the indispensable binomials of $I_G$, then
the toric ideal $I_{G}$ is generated by the indispensable
binomials. In addition the monomials of $B_i$, $1 \leq i \leq s$,
are all indispensable and also they form $T_G$.  We will prove that $$\mathcal{C}_{min} \subset
\{\supp(M_1),\ldots,\supp(M_r)\}.$$Let $E \in \mathcal{C}_{min}$ and let $\sigma=pos_{\mathbb{Q}}(A_{G})$. From Theorem 4.6 in \cite{KT} the simplicial complexes $\Delta_{G}$ and $\mathcal{D}_{\sigma}$ are identical, see \cite{KT} for the definition of the last complex. Using the fact that $I_{G}$ is generated by the binomials $B_{1},\ldots,B_{s}$ and Corollary
5.7 in \cite{KMT} we take that there is a monomial $M_{i}$ such that ${\rm cone}(M_{i})=\sigma_{E}$. For the definition and results about the cone of a monomial see \cite{KMT}. But ${\rm cone}(M_{i})=\sigma_{{\supp}(M_{i})}$, since in this case all vectors belong to an extreme ray of $\sigma$, so $\sigma_{E}=\sigma_{{\supp}(M_{i})}$ and therefore $E={\supp}(M_{i})$. Thus
$$\mathcal{C}_{min} \subset
\{\supp(M_1),\ldots,\supp(M_r)\}.$$Consider now a set
$E=\supp(M_i)$, $1 \leq i \leq r$, and we will prove that it also
belongs to $\mathcal{C}_{min}$. Suppose not, then there is an $E'
\subsetneqq E$ such that $E'=\supp({\bf x}^{{\bf u}_+})$ or
$E'=\supp({\bf x}^{{\bf u}_-})$ where ${\bf x}^{{\bf u}_+}-{\bf
x}^{{\bf u}_-} \in I_{G}$ is a circuit. Without loss of generality
we can assume that $E'=\supp({\bf x}^{{\bf u}_+})$. The monomials
${\bf x}^{{\bf u}_+}$, ${\bf x}^{{\bf u}_-}$ are squarefree and
also ${\bf x}^{{\bf u}_+}$ divides $M_i$, since $E' \subsetneqq
E$, a contradiction to the
fact that $M_i$ is indispensable.\\

Assume now that $I_{G}$ is generated by quadratic binomials. Let
$\{B_1,\ldots,B_s\}$ be a quadratic set of generators of $I_G$ and
let $S_{G}$ be the set of monomials appearing in the binomials
$B_1,\ldots,B_s$. We will prove that $S_{G}$ coincides with
$T_{G}$, i.e. $S_{G}$ is the minimal generating set of the ideal
$\mathcal{N}_{G}$. Every monomial $N$ of $S_{G}$ belongs to the
ideal $\mathcal{N}_{G}$. On the other hand for a monomial ${\bf
x}^{\bf u}\in \mathcal{N}_{G}$, there exists a monomial ${\bf
x}^{\bf w}$ such that ${\bf x}^{\bf u}-{\bf x}^{\bf w} \in I_G$.
But $I_{G}=(B_{1},\ldots,B_{s})$, so there is a monomial $N' \in
S_{G}$ which divides ${\bf x}^{\bf u}$ and therefore $S_{G}$ is a
set of generators for the ideal $\mathcal{N}_{G}$. In addition
$S_{G}$ is a minimal generating set, since every monomial $N$ of
$S_{G}$ is
quadratic. Thus $S_{G}$ is the set of indispensable monomials.\\
Using the fact that $I_G$ is generated by the binomials
$B_{1},\ldots,B_{s}$ and Corollary 5.7 in \cite{KMT} we
can easily prove that $$\mathcal{C}_{min} \subset
\{\supp(M_1),\ldots,\supp(M_r)\}.$$ It remains to prove that
$$\{\supp(M_1),\ldots,\supp(M_r)\} \subset \mathcal{C}_{min}.$$
Let $E=\supp(M_i)$, $1 \leq i \leq r$, and assume that $E$ does
not belong to $\mathcal{C}_{min}$. Then there is an $E'
\subsetneqq E$, i.e. $E'=\{i\}$ is a singleton, and a circuit
$x_{i}^{g_i}-N_i \in I_G$ such that $E'=\supp(x_{i}^{g_i})$ and
$E' \cap \supp(N_i)=\emptyset$. Let $R=A_{G}-\{{\bf a}_i\}
\subsetneqq A_{G}$ then ${\rm deg}_{G}(N_i)$ belongs to
$pos_{\mathbb{Q}}(R)$, so $g_{i}{\bf a}_i$ belongs also to
$pos_{\mathbb{Q}}(R)$, since $g_{i}{\bf a}_i={\rm deg}_{G}(N_i)$,
and therefore ${\bf a}_i \in pos_{\mathbb{Q}}(R)$. Thus
$pos_{\mathbb{Q}}(A_G)=pos_{\mathbb{Q}}(R)$, a contradiction to
the fact that $A_G$ is extremal vector configuration. \qed\\

Let $\mathcal{F} \subset I_{G}$ be a set of binomials. We shall
denote by $S({\bf b})_{\mathcal{F}}$ the graph with vertices the
elements of $\deg_{G}^{-1}({\bf b})=\{{\bf x}^{{\bf u}} \ | \
\deg_{G}({\bf x}^{{\bf u}})={\bf b}\}$ and edges the sets $\{{\bf
x}^{{\bf u}},{\bf x}^{{\bf v}}\}$ whenever ${\bf x}^{{\bf u}}-{\bf
x}^{{\bf v}}$ is a monomial multiple of a binomial in
$\mathcal{F}$. The next theorem computes the binomial arithmetical
rank and the $G$-homogeneous arithmetical rank of $I_G$, for a
bipartite graph $G$.

\begin{thm1} Let $G$ be a bipartite graph, then ${\rm
bar}(I_G)=\mu(I_G)$ and ${\rm ara}_{G}(I_G)=\mu(I_G)$.
\end{thm1}
\noindent \textbf{Proof.} First we will prove that $\{E,E'\}$ is
an edge of $\Delta_{G}$ if and only if there is an indispensable
binomial ${\bf x}^{{\bf u}_+}- {\bf x}^{{\bf u} _-} \in I_G$ with
$\supp({\bf u}_+)=E$ and $\supp({\bf u}_-)=E'$. The one
implication is easy. Let ${\bf x}^{{\bf u}_+}- {\bf x}^{{\bf u}
_-} \in I_G$ be an indispensable binomial with $\supp({\bf
u}_+)=E$ and $\supp({\bf u}_-)=E'$. Then ${\bf x}^{{\bf u}_+}-
{\bf x}^{{\bf u} _-}=f_{\Gamma}$, for an even cycle $\Gamma$ of
$G$ with no chord. But $f_{\Gamma}$ is a circuit and therefore
$\{E,E'\}$ is an edge. Conversely consider an edge $\{E,E'\}$ of
$\Delta_G$, then there is a circuit ${\bf x}^{{\bf u}_+}- {\bf
x}^{{\bf u} _-} \in I_{G}$ such that $\supp({\bf u}_+)=E$ and
$\supp({\bf u}_-)=E'$. Let $T_G=\{M_1,\ldots,M_r\}$ be the set of
indispensable monomials, then, from Lemma 3.1, there are
indispensable monomials $M_i$, $M_j$ such that $E=\supp(M_{i})$
and $E'=\supp(M_{j})$. But ${\bf x}^{{\bf u}_+}$, ${\bf x}^{{\bf
u}_-}$ are squarefree and also $M_i$, $M_j$ are minimal generators
of $\mathcal{N}_G$, so ${\bf x}^{{\bf u}_+}=M_{i}$, ${\bf x}^{{\bf
u}_-}=M_{j}$ and therefore both monomials ${\bf x}^{{\bf u}_+}$,
${\bf x}^{{\bf u}_-}$ are indispensable. Let $${\bf
b}=\deg_{G}({\bf x}^{{\bf u}_+})= \deg_{G}({\bf x}^{{\bf u}
_-}).$$ If $B_{1},\ldots,B_{s}$ are the indispensable binomials of
$I_{G}$, then $\mathcal{F}:=\{B_{1},\ldots,B_{s}\}$ is a
generating set of $I_{G}$ and therefore the graph $S({\bf
b})_{\mathcal{F}}$ is connected, see Theorem 3.2 in \cite{DS}.
Suppose that ${\bf x}^{{\bf u}_+}- {\bf x}^{{\bf u} _-}$ is not
indispensable, then there exist a path
$$(\{{\bf x}^{{\bf u}_+}={\bf x}^{{\bf u}_0},{\bf x}^{{\bf
u}_1}\},\{{\bf x}^{{\bf u}_1},{\bf x}^{{\bf u}_2}\},\ldots,\{{\bf
x}^{{\bf u}_{t-1}},{\bf x}^{{\bf u}_t}={\bf x}^{{\bf u}_-}\}), \ t
\geq 2,$$ in $S({\bf b})_{\mathcal{F}}$ connecting the vertices
${\bf x}^{{\bf u}_0}$ and ${\bf x}^{{\bf u}_t}$. Consider now the
binomial ${\bf x}^{{\bf u}_0}-{\bf x}^{{\bf u}_1}$. There is a
binomial $B_i$ and a monomial $P$ such that ${\bf x}^{{\bf
u}_0}-{\bf x}^{{\bf u}_1}=PB_{i}$ since $\{{\bf x}^{{\bf
u}_0},{\bf x}^{{\bf u}_1}\}$ is an edge of $S({\bf
b})_{\mathcal{F}}$. If $B_{i}={\bf x}^{{\bf w}_+}-{\bf x}^{{\bf
w}_-}$, then ${\bf x}^{{\bf u}_0}=P{\bf x}^{{\bf w}_+}$ and
therefore ${\bf x}^{{\bf w}_+}$ divides ${\bf x}^{{\bf u}_0}$. But
${\bf x}^{{\bf u}_0}$ is indispensable, so $P=1$ and therefore the
binomial ${\bf x}^{{\bf u}_0}-{\bf x}^{{\bf u}_1}$ is
indispensable. Thus the monomial ${\bf x}^{{\bf u}_1}$ is
indispensable. Moreover ${\bf x}^{{\bf u}_1}-{\bf x}^{{\bf u}_2}$
is indispensable, since $\{{\bf x}^{{\bf u}_1},{\bf x}^{{\bf
u}_2}\}$ is an edge of $S({\bf b})_{\mathcal{F}}$ and ${\bf
x}^{{\bf u}_1}$ is indispensable, as well as all the binomials
${\bf x}^{{\bf u}_{i-1}}-{\bf x}^{{\bf u}_{i}}$, $3 \leq i \leq
t$. Consequently there are at least two indispensable binomials
with the same $G$-degree, contradicting Theorem 3.4 in \cite{CKT}.
Recall that the $G$-degree of a binomial ${\bf x}^{{\bf u}}- {\bf
x}^{{\bf v}} \in I_G$ is defined to be $\deg_G({\bf x}^{{\bf u}}-
{\bf x}^{{\bf v}}):=\deg_{G}({\bf
x}^{{\bf u}})$.\\
Now, from Theorem 3.4 in \cite{CKT}, every edge of $\Delta_G$
constitute a connected component. Remark that $\Delta_{G}$ has no
connected components which are singletons, since
$\mathcal{C}_{min}=\{\supp(M_1),\ldots,\supp(M_r)\}$ and $I_G$ is
generated by the indispensable binomials. Thus every connected
component of $\Delta_G$ is an edge, so $\Delta_G$ has $s$
connected components. Let $\Delta _{G}^{i}=\{E,E'\}$, $1 \leq i
\leq s$, be a connected component of $\Delta_{G}$. There are two
maximal $\{0,1\}$-matchings in $\Delta _{G}^{i}$, namely
$\{\{E,E'\}\}$ and $\{\{E\},\{E'\}\}$. We have that $\delta(\Delta
_{G}^{i})_{\{0,1\}}=1$ which is attained by the maximal
$\{0,1\}$-matching $\{\{E,E'\}\}$. Consequently
$$\delta(\Delta_{G})_{\{0,1\}}=\sum_{i=1}^{s}\delta(\Delta _{G}^{i})_{\{0,1\}}=s,$$
i.e. $\delta(\Delta_{G})_{\{0,1\}}=\mu(I_G)$. From Theorem 2.4 we have that ${\rm bar}(I_G) \geq \mu(I_G)$ and therefore ${\rm bar}(I_G) =\mu(I_G)$.\\
In addition
$\delta(\Delta_{G})_{\Omega}=\delta(\Delta_{G})_{\{0,1\}}$, since
${\rm dim}(\Delta _{G})=1$. So
$\delta(\Delta_{G})_{\Omega}=\mu(I_G)$ and therefore, from Theorem
2.4, we have that ${\rm ara}_{G}(I_G) \geq \mu(I_G)$. Thus ${\rm
ara}_{G}(I_G)=\mu(I_G)$. \qed

\begin{ex1}{\rm Let $\mathcal{K}_{3,3}$ be the complete bipartite graph on the vertex set $\{v_{1},\ldots,v_{6}\}$ with $9$ edges: $$\{v_{1},v_{4}\},\{v_{1},v_{5}\},\{v_{1},v_{6}\},\{v_{2},v_{4}\},\{v_{2},v_{5}\},\{v_{2},v_{6}\},\{v_{3},v_{4}\},\{v_{3},v_{5}\},\{v_{3},v_{6}\}.$$ We consider one variable $x_{ij}$, $1 \leq i<j \leq 6$, for each edge
$\{v_{i},v_{j}\}$ of $G$ and form the polynomial ring $K[x_{ij}| 1
\leq i<j \leq 6]$. The toric ideal $I_G$ is minimally generated by
$9$ binomials: $$x_{14}x_{26}-x_{16}x_{24},
x_{15}x_{36}-x_{16}x_{35},x_{25}x_{36}-x_{26}x_{35},x_{24}x_{36}-x_{26}x_{34},x_{14}x_{25}-x_{15}x_{24},$$
$$x_{15}x_{26}-x_{16}x_{25}, x_{24}x_{35}-x_{25}x_{34},
x_{14}x_{36}-x_{16}x_{34},x_{14}x_{35}-x_{15}x_{34}.$$ The
simplicial complex $\Delta_{G}$ has $18$ vertices, corresponding
to the $18$ monomials arising in the above minimal generating set
of $I_G$, and $9$ edges corresponding to the $9$ minimal
generators of $I_G$. From Theorem 3.2 we have that ${\rm
ara}_{G}(I_G)={\rm bar}(I_G)=9$. The height of $I_{G}$ equals
$9-6+1=4$, see Proposition 3.2 in \cite{Vil}. For the arithmetical
rank of $I_G$ we have that $4 \leq {\rm ara}(I_G) \leq 7$, since
$I_G$ equals the radical of the ideal generated by
$$x_{14}x_{26}-x_{16}x_{24}+
x_{15}x_{36}-x_{16}x_{35},x_{25}x_{36}-x_{26}x_{35}+x_{14}x_{25}-x_{15}x_{24},x_{24}x_{36}-x_{26}x_{34},$$ $$x_{15}x_{26}-x_{16}x_{25}, x_{24}x_{35}-x_{25}x_{34},
x_{14}x_{36}-x_{16}x_{34},x_{14}x_{35}-x_{15}x_{34}.$$}
\end{ex1}

\par An interesting case occurs when $I_G$ is generated by quadratic
binomials. In \cite{OH} a combinatorial criterion for the toric
ideal $I_G$ to be generated by quadratic binomials is studied.
Remark that if $B=x_{i}x_{j}-x_{k}x_{l}$ is a quadratic binomial
in $I_G$, then $B=f_{\Gamma}$ for an even cycle $\Gamma$ of $G$ of
length $4$. We are going to compute the binomial arithmetical rank
and the $G$-homogeneous arithmetical rank of such an ideal.\\
We shall denote by $\Delta_{{\rm ind}(A_{G})}$ the indispensable
complex of $A_G$. For the definition and results about the
indispensable complex see \cite{CKT}.
\begin{prop1} Let $I_G$ be a toric ideal
generated by quadratic binomials, then
\begin{enumerate} \item $\{E,E'\}$ is a connected component of
$\Delta_{G}$ if and only if there is an indispensable quadratic
binomial ${\bf x}^{{\bf u}_+}-{\bf x}^{{\bf u}_-} \in I_G$ with $\supp({\bf x}^{{\bf u}_+})=E$ and $\supp({\bf x}^{{\bf u}_-})=E'$. \item every connected component of $\Delta_{G}$
is either an edge or a $2$-simplex .
\end{enumerate}
\end{prop1}
\noindent \textbf{Proof.} (1) The first goal is to prove that
$\{E,E'\}$ is an edge of $\Delta_{G}$ if and only if there is a
quadratic binomial ${\bf x}^{{\bf u}_+}- {\bf x}^{{\bf u}_-} \in
I_G$ with $\supp({\bf x}^{{\bf u}_+})=E$ and $\supp({\bf x}^{{\bf
u}_-})=E'$. The one implication follows from the fact that if
${\bf x}^{{\bf u}_+}- {\bf x}^{{\bf u}_-} \in I_G$ is a quadratic
binomial, then $\supp({\bf u})$ is minimal with respect to
inclusion, since $A_G$ is extremal, and therefore the binomial
${\bf x}^{{\bf u}_+}- {\bf x}^{{\bf u}_-}$ is a circuit.
Conversely let $\{E,E'\}$ be an edge of $\Delta_{G}$ and let ${\bf
x}^{{\bf u}_+}- {\bf x}^{{\bf u}_-} \in I_{G}$ be a circuit with
$\supp({\bf x}^{{\bf u}_+})=E$ and $\supp({\bf x}^{{\bf
u}_-})=E'$, where ${\bf u}=(u_{1},\ldots,u_{m}) \in \mathbb{Z}^m$.
From Corollary 8.4.16 in \cite{V} we have that $|u_i| \leq 2$. Let
$\{B_1,\ldots,B_s\}$ be a quadratic set of generators of $I_G$ and
let $T_{G}=\{M_1,\ldots,M_r\}$ be the indispensable monomials. We
have that both $E$, $E'$ belong to $\mathcal{C}_{min}$ and
therefore, from Lemma 3.1, $E=\supp(M_i)$ and $E'=\supp(M_j)$. But
$M_i$ and $M_{j}$ are quadratic monomials, so $E$ and $E'$ consist
of exactly 2 elements. Let $E=\{k,l\}$ and $E'=\{p,q\}$. We will
consider three cases.
\begin{enumerate} \item[(i)] If the monomials ${\bf x}^{{\bf u}_+}$ and
${\bf x}^{{\bf u}_-}$ are squarefree, then ${\bf x}^{{\bf u}_+}-
{\bf x}^{{\bf u}_-}$ is a quadratic binomial with $\supp({\bf
x}^{{\bf u}_+})=E$ and $\supp({\bf x}^{{\bf u}_-})=E'$.
\item[(ii)] If ${\bf x}^{{\bf u}_+}=x_{k}^{2}x_{l}^{2}$ and ${\bf
x}^{{\bf u}_-}=x_{p}^{2}x_{q}^{2}$, then ${\rm deg}_{G}({\bf
x}^{{\bf u}_+})={\rm deg}_{G}({\bf x}^{{\bf u}_-})$ and therefore
$2{\bf a}_k+2{\bf a}_l=2{\bf a}_p+2{\bf a}_q$. Thus ${\bf
a}_k+{\bf a}_l={\bf a}_p+{\bf a}_q$, so the binomial
$x_{k}x_{l}-x_{p}x_{q}$ belongs to $I_{G}$ and also
$\supp(x_{k}x_{l})=E$, $\supp(x_{p}x_{q})=E'$. \item[(iii)] If
${\bf x}^{{\bf u}_+}=x_{k}^{2}x_{l}$ and ${\bf x}^{{\bf
u}_-}=x_{p}^{2}x_{q}$, then ${\rm deg}_{G}({\bf x}^{{\bf
u}_+})={\rm deg}_{G}({\bf x}^{{\bf u}_-})$ and therefore $2{\bf
a}_k+{\bf a}_l=2{\bf a}_p+{\bf a}_q$. Assume that ${\bf a}_k-{\bf
a}_p \neq {\bf 0}$. Notice that every nonzero coordinate of the
above vector equals either $1$ or $-1$. We have that ${\bf
a}_q-{\bf a}_l=2({\bf a}_k-{\bf a}_p)$, so every nonzero
coordinate of the vector ${\bf a}_q-{\bf a}_l$ equals either $2$
or $-2$, contradiction. Thus ${\bf a}_k-{\bf a}_p={\bf 0}$, so
${\bf a}_{k}={\bf a}_{p}$ and therefore $k=p$, a contradiction to
the fact that $\supp({\bf u}_+) \cap \supp({\bf
u}_-)=\emptyset$.\\Similarly the assumption ${\bf x}^{{\bf
u}_+}=x_{k}x_{l}^{2}$ or ${\bf x}^{{\bf u}_-}=x_{p}x_{q}^{2}$
leads again to a contradiction.
\end{enumerate}
The second goal is to prove that $\{E,E'\}$ is a connected
component of $\Delta_{G}$ if and only if there is an indispensable
quadratic binomial $M_{i}- M_{j} \in I_G$ with $\supp(M_i)=E$ and
$\supp(M_j)=E'$. Suppose that the binomial $M_{i}- M_{j} \in
I_{G}$ is indispensable with $\supp(M_i)=E$ and $\supp(M_j)=E'$.
From Theorem 3.4 in \cite{CKT} we have that $\{M_{i},M_{j}\}$ is a
facet of the indispensable complex $\Delta_{{\rm ind}(A_{G})}$.
Assume that $\{E,E'\}$ is not a connected component of
$\Delta_{G}$. Let us suppose that $\{E,E''\}$ is an edge of
$\Delta_{G}$, then there exists a quadratic binomial ${\bf
x}^{{\bf u}_+}- {\bf x}^{{\bf u}_-} \in I_G$ with $\supp({\bf
x}^{{\bf u}_+})=E$ and $\supp({\bf x}^{{\bf u}_-})=E''$. Moreover
${\bf x}^{{\bf u}_+}=M_{i}$ and ${\bf x}^{{\bf u}_-}=M_{k}$, since
the monomials ${\bf x}^{{\bf u}_+}$, ${\bf x}^{{\bf u}_-}$ are
quadratic and therefore indispensable. Notice that ${\rm
deg}_{G}(M_{i})={\rm deg}_{G}(M_{k})$. Thus $${\rm
deg}_{G}(M_{i})={\rm deg}_{G}(M_{j})={\rm deg}_{G}(M_{k}),$$ since
the binomial $M_{i}- M_{j}$ belongs to $I_G$ and therefore ${\rm
deg}_{G}(M_{i})={\rm deg}_{G}(M_{j})$. So $\{M_{i},M_{j},M_{k}\}$
is a face of $\Delta_{{\rm ind}(A_{G})}$, a contradiction to the
fact that $\{M_{i},M_{j}\}$ is a facet of $\Delta_{{\rm
ind}(A_{G})}$. Consequently $\{E,E'\}$ is a connected component of
$\Delta_{G}$. Conversely assume that $\{E,E'\}$ is a connected component of $\Delta_{G}$. Then there is a quadratic binomial ${\bf x}^{{\bf u}_+}- {\bf x}^{{\bf u}_-} \in I_G$ with $\supp({\bf x}^{{\bf u}_+})=E$ and $\supp({\bf x}^{{\bf u}_-})=E'$. In fact ${\bf x}^{{\bf u}_+}=M_{i}$ and ${\bf x}^{{\bf u}_-}=M_{j}$, for some indispensable monomials $M_{i}$, $M_{j}$. Suppose that the above binomial is not indispensable, then, since $M_{i}$ is indispensable, there is an $l \in \{1,\ldots,s\}$ such that $B_{l}=M_{i}-M_{k}$. Set $E''=\supp(M_{k}) \in \mathcal{C}_{min}$. We have that $\{E,E''\}$ and $\{E',E''\}$ are edges of $\Delta_{G}$, since also $M_{j}-M_{k} \in I_{G}$, a contradiction to the fact that $\{E,E'\}$ is connected component.\\
(2) Notice that $\Delta_{G}$ has no connected components which are
singletons. To see this consider $E=\supp(M_i) \in
\mathcal{C}_{min}$, then there is an $l \in \{1,\ldots,s\}$ such
that $B_{l}=M_{i}-M_{j}$. Consequently $\{\supp(M_i),\supp(M_j)\}$
is an edge of $\Delta_{G}$.\\ Next we will show that
$\{E,E',E''\}$ is a $2$-simplex of $\Delta_{G}$ if and only if
there are quadratic binomials $M_{i}-M_{j}$, $M_{j}-M_{k}$,
$M_{i}-M_{k}$ in $I_G$ with $\supp(M_{i})=E$, $\supp(M_{j})=E'$
and $\supp(M_{k})=E''$. The one implication is easily derived from
the fact that if $\{E,E',E''\}$ is a $2$-simplex of $\Delta_{G}$,
then every $2$-element subset of it is an edge. Conversely we have
that
$${\rm deg}_{G}(M_i)={\rm deg}_{G}(M_j)={\rm deg}_{G}(M_k)$$ belongs to the
intersection $$relint_{\mathbb{Q}}\left( \sigma_{E}\right) \cap
relint_{\mathbb{Q}} \left( \sigma_{E'}\right) \cap
relint_{\mathbb{Q}}\left( \sigma_{E''}\right).$$Thus
$\{E,E',E''\}$ is a $2$-simplex of $\Delta_{G}$.\\Finally we prove
that if $\{E,E',E''\}$ is a $2$-simplex of $\Delta_{G}$, then it
is a connected component. There are quadratic binomials
$M_{i}-M_{j}$, $M_{j}-M_{k}$ and $M_{i}-M_{k}$ in $I_G$, where
$\supp(M_i)=E$, $\supp(M_j)=E'$ and $\supp(M_k)=E''$. Let us
suppose that $M_{i}=x_{i_1}x_{i_2}$, $M_{j}=x_{j_1}x_{j_2}$ and
$M_{k}=x_{k_1}x_{k_2}$. In addition there are even cycles
$\Gamma_{1}$, $\Gamma_{2}$, $\Gamma_{3}$ of $G$ of length $4$ such
that $M_{i}-M_{j}=f_{\Gamma_{1}}$, $M_{i}-M_{k}=f_{\Gamma_{2}}$
and $M_{j}-M_{k}=f_{\Gamma_{3}}$. The cycle $\Gamma_{1}$ has $4$
edges, namely $e_{i_{1}}$, $e_{i_{2}}$, $e_{j_{1}}$ and
$e_{j_{2}}$, the cycle $\Gamma_{2}$ has $4$ edges, namely
$e_{i_{1}}$, $e_{i_{2}}$, $e_{k_{1}}$ and $e_{k_{2}}$, and
$\Gamma_{3}$ has $4$ edges, namely $e_{j_{1}}$, $e_{j_{2}}$,
$e_{k_{1}}$ and $e_{k_{2}}$. Notice that the edges $e_{i_{1}}$ and
$e_{i_{2}}$ have no common vertex. The above three cycles have the
same vertex set $\mathcal{V}$ consisting of $4$ vertices. Moreover
these are the only even cycles of length $4$ with vertex set
$\mathcal{V}$. Let $\mathcal{K}_4$ be the induced subgraph of $G$
on the above vertex set. It is a complete subgraph with $4$
vertices and edges
$$E(\mathcal{K}_4)=\{e_{i_1},e_{i_2},e_{j_1},e_{j_2},e_{k_1},e_{k_2}\}.$$ If, for
example, $\{E,E'''\}$ is an edge of $\Delta_{G}$, then there is a
quadratic binomial $M_{i}-M_{l} \in I_G$ with $\supp(M_{l})=E'''
\in C_{min}$. Furthermore $M_{i}-M_{l}=f_{\Gamma_{4}}$, for an
even cycle $\Gamma_{4}$ of $G$ of length $4$. The vertex set of
$\Gamma_{4}$ is $\mathcal{V}$, since $e_{i_{1}}$ and
$e_{i_{2}}$ have no common vertex, and therefore $\Gamma_{4}$ coincides with either $\Gamma_{1}$ or $\Gamma_{2}$. Thus $M_{l}$ equals to either $M_{j}$ or $M_{k}$, so $E'''=E'$ or $E'''=E''$. Consequently $\{E,E',E''\}$ is a connected component of $\Delta_{G}$. \qed \\

\begin{rem1} {\rm (1) To every
connected component $\{E,E'\}$ of $\Delta_{G}$ we can associate an indispensable binomial $f_{\Gamma}=x_{i}x_{j}-x_{k}x_{l} \in I_G$, for an even
cycle $\Gamma$ of $G$ of length $4$, where $E=\{i,j\}$ and $E'=\{k,l\}$, and also the induced subgraph $\mathcal{H}$ of $G$ on the vertex set of $\Gamma$. The subgraph $\mathcal{H}$ is not a complete graph. Moreover the toric $I_{\mathcal{H}}$ is complete intersection of height $1$ and it is generated by $f_{\Gamma}$.\\
(2) $\{E,E',E''\}$ is a $2$-simplex of $\Delta_{G}$ if and only if
there are quadratic binomials $M_{i}-M_{j}$, $M_{i}-M_{k}$,
$M_{j}-M_{k}$ in $I_G$
with $\supp(M_{i})=E$, $\supp(M_{j})=E'$ and $\supp(M_{k})=E''$.\\
(3) To every connected component of $\Delta_{G}$ which is a
$2$-simplex we can associate a complete subgraph $\mathcal{K}_4$
of $G$ of order $4$. The toric ideal $I_{\mathcal{K}_4}$ is
minimally generated by two binomials $f_{\Gamma_{1}}$ and
$f_{\Gamma_{2}}$, where $\Gamma_{1}$ and $\Gamma_{2}$ are even
cycles of length $4$ on the vertex set of $\mathcal{K}_4$.}
\end{rem1}
\begin{prop1} Let $\Gamma=(e_{i},e_{p},e_{j},e_{q})$ be an even
cycle of a graph $G$ such that the induced subgraph $\mathcal{H}$
of $G$ on the vertex set of $\Gamma$ is not a complete graph. If
$H$ is a nonzero polynomial in $I_{\mathcal{H}}$, then there exist
monomials $M$, $N$ of $H$ such that $x_{i}x_{j}$ divides $M$ and
$x_{p}x_{q}$ divides $N$.
\end{prop1}
\noindent \textbf{Proof.} For the toric ideal $I_{\mathcal{H}}$ we
have, from Proposition 4.13 in \cite{St}, that
$I_{\mathcal{H}}=I_{G} \cap K[x_l | e_{l} \in E(\mathcal{H})]$.
Since $I_{\mathcal{H}}=(f_{\Gamma})$, there is a nonzero
polynomial $C \in K[x_l | e_{l} \in E(\mathcal{H})]$ such that
$H=Cf_{\Gamma}$. The polynomial $C$ has a unique representation as
a sum of terms $C=C_{1}+\cdots+C_{s}$. Notice that the monomials
of two distinguished terms $C_k$ and $C_{l}$ are different. We
have that
\begin{equation*} \tag{3.1}
  H=C_{1}x_{i}x_{j}+\cdots+C_{s}x_{i}x_{j}-C_{1}x_{p}x_{q}-\cdots-C_{s}x_{p}x_{q}.
 \end{equation*}
Assume that $H$ has no term whose monomial is $x_{p}x_{q}$. This
implies that in the above expression of $H$ all the terms of the
form $C_{k}x_{p}x_{q}$ should by cancelled. But these terms can
not cancel by themselves, so terms of the form $C_{k}x_{i}x_{j}$
are used to cancel them. We claim that a term $C_{k}x_{i}x_{j}$
can be used to cancel atmost one term of the form
$-C_{l}x_{p}x_{q}$. Assume that $C_{k}x_{i}x_{j}$ is used to
cancel the terms $-C_{l}x_{p}x_{q}$ and $-C_{r}x_{p}x_{q}$. Let
$M_{1}$, $M_{2}$ and $M_{3}$ be the monomials of the terms
$C_{k}x_{i}x_{j}$, $-C_{l}x_{p}x_{q}$ and $-C_{r}x_{p}x_{q}$
respectively. Then $M_{1}=M_{2}$ and $M_{1}=M_{3}$, so
$M_{2}=M_{3}$, contradiction. But $x_{p}x_{q}$ divides no
monomials of $H$, so there are two cases.
\begin{enumerate} \item In expression (3.1) every term cancels.
Therefore $H$ is equal to zero, contradiction. \item In expression
(3.1) every term of the form $-C_{l}x_{p}x_{q}$ cancels, but still
there
exist terms of the form $C_{k}' x_{i}x_{j}$ where $C_{k}'$ is different from $C_k$. Notice that the monomial of such a term coincides with the monomial of a suitable term $-C_{l}x_{p}x_{q}$. Thus every term $C_{k}' x_{i}x_{j}$ is divided by $x_{p}x_{q}$, contradiction. \qed\\
\end{enumerate}

The next lemma will be useful in the proof of Theorem 3.8.

\begin{lem1} Let $G$ be a graph with edges $E(G)=\{e_{1},\ldots,e_{m}\}$, $\Gamma$ an even cycle of length
$4$ and $\mathcal{H}$ the induced subgraph of $G$ on the vertex set of $\Gamma$.
If $\mathcal{F} \subset I_G$ is a set of
$G$-homogeneous polynomials which generates $I_G$ up to radical,
then $\mathcal{F} \cap K[x_{i}|e_{i} \ \textrm{is an edge of} \
\mathcal{H}]$ generates $I_{\mathcal{H}}$ up to radical.
\end{lem1}
\noindent \textbf{Proof.} Let $\{v_{i_{1}},\ldots,v_{i_{4}}\}$ be the vertices of $\Gamma$. The rational polyhedral cone
$pos_{\mathbb{Q}}(A_{\mathcal{H}})$ is
a face of $pos_{\mathbb{Q}}(A_G)$ with defining vector ${\bf c}=(c_{1},\ldots,c_{n}) \in \mathbb{Z}^n$ having coordinates $$c_j= \left \{\begin{array} {lll} 0, & \textrm{if} \ j=i_{1},i_{2},i_{3},i_{4}\\
 1, & \textrm{otherwise}. \end{array} \right.$$Thus, from
 Proposition 3.2 in \cite{KMT1}, we have that $\mathcal{F} \cap K[x_{i}|e_{i}
 \in
E(\mathcal{H})]$ generates $I_{\mathcal{H}}$ up to radical. \qed\\

The following theorem determines the binomial arithmetical rank
and the $G$-homogeneous arithmetical rank of a toric ideal $I_G$
generated by quadratic binomials.
\begin{thm1} Let $G$ be a graph. If $I_G$ is generated by quadratic
binomials, then \begin{enumerate} \item ${\rm bar}(I_G)=\mu(I_G)$
and \item ${\rm ara}_{G}(I_G)=\mu(I_G)$. \end{enumerate}
\end{thm1}
\noindent \textbf{Proof.} (1) Let $g \geq 0$ be the number of
indispensable binomials of $I_G$, then, from Proposition 3.4 (1),
the simplicial complex $\Delta_{G}$ has exactly $g$ connected
components which are edges.\\ We will show that $\Delta_{G}$ has
$\frac{s-g}{2}$ connected components which are $2$-simplices,
where $s=\mu(I_G)$. Let $\mathcal{B}=\{B_{1},\ldots,B_{s}\}$ be a
minimal set of quadratic generators of $I_{G}$ and let
$\{M_{1},\ldots,M_{r}\}$ be the set of indispensable monomials.
Notice that $\mathcal{B}$ has $s-g$ binomials which are not
indispensable. Given a $2$-simplex $\{E,E',E''\}$ of $\Delta_{G}$,
there are quadratic binomials $f_{\Gamma_{1}}=M_{i}-M_{j}$,
$f_{\Gamma_{2}}=M_{i}-M_{k}$, $f_{\Gamma_{3}}=M_{j}-M_{k}$ in
$I_G$ with $\supp(M_{i})=E$, $\supp(M_{j})=E'$ and
$\supp(M_{k})=E''$. Remark that the binomials $f_{\Gamma_{1}}$,
$f_{\Gamma_{2}}$ and $f_{\Gamma_{3}}$, as well as
$-f_{\Gamma_{1}}$, $-f_{\Gamma_{2}}$ and $-f_{\Gamma_{3}}$, are
not indispensable. In fact the minimal generating set
$\mathcal{B}$ contains exactly two of the binomials
$f_{\Gamma_{1}}$, $-f_{\Gamma_{1}}$, $f_{\Gamma_{2}}$,
$-f_{\Gamma_{2}}$, $f_{\Gamma_{3}}$ and $-f_{\Gamma_{3}}$, since
the monomials $M_{i}$, $M_{j}$, $M_{k}$ are indispensable and
$\{E,E',E''\}$ is a connected component of $\Delta_{G}$. Let $t$
be the number of connected components which are $2$-simplices,
then $\mathcal{B}$ contains at least $2t$ binomials which are not
indispensable. So $2t \leq s-g$. On the other hand if
$B_{l}=M_{i}-M_{j}$ is not indispensable, then
$\{\supp(M_i),\supp(M_{j})\}$ is an edge which is not a connected
component of $\Delta_{G}$. Therefore there is a monomial $M_{k}$
such that $\{\supp(M_i),\supp(M_{j}),\supp(M_{k})\}$ is a
connected component of $\Delta_{G}$. Moreover there exists a $p
\in \{1,\ldots,s\}$ such that $B_p$ or $-B_p$ equals either
$M_{i}-M_{k}$ or $M_{j}-M_{k}$. But $\mathcal{B}$ is a minimal
generating set, so there exist exactly two binomials in
$\mathcal{B}$ whose monomials are $M_{i}$, $M_{j}$ and $M_{k}$.
Thus $\Delta_{G}$ has at least $\frac{s-g}{2}$ connected
components which are $2$-simplices, so $\frac{s-g}{2} \leq t$.
Consequently $t=\frac{s-g}{2}$.\\ For every connected component
$\Delta_{G}^{i}$ of $\Delta_{G}$ which is an edge we have
$\delta(\Delta_{G}^{i})_{\{0,1\}}=1$, while for every connected
component $\Delta_{G}^{i}$ of $\Delta_{G}$ which is a $2$-simplex
we have $\delta(\Delta_{G}^{i})_{\{0,1\}}=2$. Consequently
$$\delta(\Delta_{G})_{\{0,1\}}=g+2\frac{s-g}{2}=s,$$ i.e. $\delta(\Delta_{G})_{\{0,1\}}=\mu(I_G)$, and therefore, from Theorem 2.4, ${\rm bar}(I_G)=\mu(I_G)$.\\
(2) Let $\mathcal{F} \subset I_G$ be a set of $G$-homogeneous
polynomials which generate $I_G$ up to radical. Let
$\Delta_{G}^{i}=\{E,E'\}$, $\Delta_{G}^{j}$ be two connected
components which are edges and let $\mathcal{H}_i$ and
$\mathcal{H}_j$, respectively, be the corresponding induced
subgraphs. Let $E=\{k,l\}$ and $E'=\{p,q\}$, then
$I_{\mathcal{H}_i}=(f_{\Gamma})$ where
$f_{\Gamma}=x_{k}x_{l}-x_{p}x_{q}$. The cycle $\Gamma$ has $4$
edges, namely $e_{k}$, $e_{l}$, $e_{p}$ and $e_{q}$. From
Proposition 3.6 every nonzero $H \in I_{\mathcal{H}_i}$ is a
polynomial in at least 4 variables, namely $x_{k}$, $x_{l}$,
$x_{p}$ and $x_{q}$. We will prove that every nonzero polynomial
$H$, which belongs to $I_{\mathcal{H}_i}$, does not belong to
$I_{\mathcal{H}_j}$. Assume that there is a nonzero polynomial $H
\in I_{\mathcal{H}_i}$ which belongs to $I_{\mathcal{H}_j}$. From
Proposition 4.13 in \cite{St}, we have that
$I_{\mathcal{H}_j}=I_{G} \cap K[x_{i}|e_{i} \in
E(\mathcal{H}_j)]$. But $H$ belongs to $I_{\mathcal{H}_j}$, so $H$
is a polynomial in the ring $K[x_{i}|e_{i} \in E(\mathcal{H}_j)]$
and therefore every edge of $\Gamma$ is also an edge of
$\mathcal{H}_j$. Thus the indispensable binomial $f_{\Gamma}$
belongs to $I_{\mathcal{H}_j}$ and therefore, from Proposition 3.4
(1), we have that $\{E,E'\}$ is a connected component of
$\Delta_{G}^{j}$, a contradiction. Given a connected component of
$\Delta_{G}$, which is an edge, and the corresponding induced
subgraph $\mathcal{H}$ of $G$, there exists, from Lemma 3.7, at
least one $G$-homogeneous polynomial $F \in I_{\mathcal{H}}$ in
$\mathcal{F}$. The simplicial complex $\Delta_{G}$ has $g$
connected components which are edges, so $\mathcal{F}$ has at
least $g$ $G$-homogeneous polynomials, say $F_{1},\ldots,F_{g}$,
belonging to the corresponding toric ideals $I_{\mathcal{H}_i}$,
$1 \leq i \leq g$.\\ Remark that if $G$ has a complete subgraph
$\mathcal{K}_4$, then every polynomial $F_{i}$, $1 \leq i \leq g$,
does not belong to the toric ideal $I_{\mathcal{K}_4}$, since
$I_{\mathcal{H}_i}=I_{G} \cap K[x_{r}|e_{r} \in
E(\mathcal{H}_i)]$ and every $\mathcal{H}_i$ is not a complete graph.\\
Let $\Delta_{G}^{p}$, $\Delta_{G}^{q}$ be two connected components
which are $2$-simplices and let $\mathcal{K}_{4,p}$ and
$\mathcal{K}_{4,q}$, respectively, the corresponding induced
subgraphs. We will prove that every nonzero polynomial $H$, which
is in the ideal $I_{\mathcal{K}_{4,p}}$, does not belong to
$I_{\mathcal{K}_{4,q}}$. Let
$I_{\mathcal{K}_{4,p}}=(f_{\Gamma_{1}},f_{\Gamma_{2}})$, where
$f_{\Gamma_{1}}=x_{i_1}x_{i_2}-x_{i_3}x_{i_4}$ and
$f_{\Gamma_{2}}=x_{i_1}x_{i_2}-x_{i_5}x_{i_6}$. Assume that $H \in
I_{\mathcal{K}_{4,p}}$ is a nonzero polynomial which belongs to
$I_{\mathcal{K}_{4,q}}$. Since $H$ belongs to
$I_{\mathcal{K}_{4,p}}$, every monomial of $H$ is of the form
$Cx_{_{i_1}x_{i_2}}$ or $Nx_{_{i_3}x_{i_4}}$ or
$Qx_{_{i_5}x_{i_6}}$, for appropriate monomials $C$, $N$ and $Q$.
But $I_{\mathcal{K}_{4,q}}=I_{G} \cap K[x_{j}|e_{j} \in
E(\mathcal{K}_{4,q})]$, so $\mathcal{K}_{4,q}$ has at least $2$
edges coming from $\mathcal{K}_{4,p}$. These edges are $e_{i_1}$
and $e_{i_2}$ or $e_{i_3}$ and $e_{i_4}$ or $e_{i_5}$ and
$e_{i_6}$. Notice that the edges $e_{i_1}$ and $e_{i_2}$ have no
common vertex. The same holds for $e_{i_3}$ and $e_{i_4}$, as well
as the edges $e_{i_5}$ and $e_{i_6}$. But $\mathcal{K}_{4,q}$ is a
complete graph, so $\mathcal{K}_{4,p}=\mathcal{K}_{4,q}$
contradiction.\\Given a connected component of $\Delta_{G}$, which
is a $2$-simplex, and the corresponding induced subgraph
$\mathcal{K}_{4}$ of $G$, there exist, from Lemma 3.7, at least
two $G$-homogeneous polynomial $H_{1},H_{2} \in
I_{\mathcal{K}_{4}}$ in $\mathcal{F}$. The simplicial complex
$\Delta_{G}$ has $\frac{s-g}{2}$ connected components which are
$2$-simplices, so $\mathcal{F}$ has also at least
$2\frac{s-g}{2}=s-g$ $G$-homogeneous polynomials, say
$H_{1},\ldots,H_{s-g}$, belonging to the corresponding toric
ideals $I_{\mathcal{K}_{4,i}}$, $1 \leq i \leq s-g$. Thus $${\rm
ara}_{G}(I_G) \geq g+(s-g)=s,$$ i.e. ${\rm ara}_{G}(I_G) \geq
\mu(I_G)$, and therefore ${\rm ara}_{G}(I_G)=\mu(I_G)$. \qed

\begin{ex1} {\rm Consider the complete graph $\mathcal{K}_n$, $n \geq 4$, on the vertex set $\{v_{1},\ldots,v_{n}\}$. We consider one variable $x_{ij}$, $1 \leq i<j \leq n$, for each edge
$\{v_{i},v_{j}\}$ of $\mathcal{K}_n$ and form the polynomial ring
$K[x_{ij}| 1 \leq i<j \leq n]$. The toric ideal
$I_{\mathcal{K}_n}$ is the kernel of the $K$-algebra homomorphism
$$ \phi: K[x_{ij}| 1 \leq i<j \leq n]\rightarrow
K[t_1,\dots,t_{n}]$$ given by
$$\phi(x_{ij}) = t_it_j.$$From Proposition 3.2 in \cite{Vil} the height of $I_{\mathcal{K}_n}$ equals ${{n}\choose {2}}-n=\frac{n(n-3)}{2}$, i.e. the number of edges minus the number of vertices. It is well known, see for example Proposition 9.2.1 in \cite{V}, that
$$B=\{x_{ij}x_{kl}-x_{il}x_{jk}, x_{ik}x_{jl}-x_{il}x_{jk} | 1
\leq i<j<k<l \leq n\}$$ is a minimal generating set for
$I_{\mathcal{K}_n}$. The toric ideal $I_{\mathcal{K}_n}$ has no
indispensable binomials and therefore every connected component of $\Delta_{\mathcal{K}_n}$ is a $2$-simplex. Thus $\Delta_{\mathcal{K}_n}$ has
$3{{n}\choose {4}}$ vertices and ${n}\choose {4}$ connected
components, which are $2$-simplices, corresponding to all complete
subgraphs of $\mathcal{K}_n$ of order $4$. For the minimal number
of generators we have that
$$\mu(I_{\mathcal{K}_n})=2 {{n}\choose {4}}=\frac{n(n-1)(n-2)(n-3)}{12}.$$
Consequently $${\rm bar}(I_{\mathcal{K}_n})={\rm
ara}_{\mathcal{K}_n}(I_{\mathcal{K}_n})=\frac{n(n-1)(n-2)(n-3)}{12}.$$
Using the result of Eisenbud-Evans and Storch that ${\rm
ara}(I_{\mathcal{K}_n})$ is bounded above by the number of
variables of $K[x_{ij}| 1 \leq i<j \leq n]$ we take that $$
\frac{n(n-3)}{2} \leq {\rm ara}(I_{\mathcal{K}_n}) \leq
\frac{n(n-1)}{2}.$$For the polynomials which minimally generate
$I_{\mathcal{K}_n}$ up to radical we know, from Theorem 5.8 in
\cite{KMT}, that there must be at least $3{{n}\choose {4}}$
monomials in at least $2{{n}\choose {4}}$
$\mathcal{K}_n$-homogeneous components.\\}
\end{ex1}

\noindent {\bf Acknowledgment.} The author would like to thank the
referee for helpful comments.

\end{document}